\documentclass[a4paper,11pt,leqno]{amsart}

\usepackage{amsmath,amsthm,amssymb}
\usepackage{graphicx,xcolor}

\numberwithin{equation}{section}
\theoremstyle{plain}
\newtheorem{theorem}{Theorem}[section]
\newtheorem{introtheorem}{Theorem}

\newtheorem{proposition}[theorem]{Proposition}
\newtheorem{fact}[theorem]{Fact}
\newtheorem{lemma}[theorem]{Lemma}
\newtheorem{corollary}[theorem]{Corollary}

\theoremstyle{remark}
\newtheorem{definition}[theorem]{Definition}
\newtheorem{remark}[theorem]{Remark}
\newtheorem*{remark*}{Remark}
\newtheorem{example}[theorem]{Example}
\newtheorem*{acknowledgements}{Acknowledgements}

\newcommand{\C}{{\mathbb C}}
\newcommand{\R}{\mathbb{R}}
\newcommand{\mc}[1]{\mathcal{#1}}
\newcommand{\mb}[1]{\mathbf{#1}}
\newcommand{\op}{\operatorname}
\newcommand{\pmt}[1]{{\begin{pmatrix}#1\end{pmatrix}}}
\newcommand{\pt}[1]{\text{\rm #1}}
\renewcommand{\phi}{\varphi}
\renewcommand{\epsilon}{\varepsilon}

\title[Entire flat Chebyshev metrics]{
Entire flat Chebyshev metrics
and applications to zero mean curvature surfaces
in Lorentz-Minkowski space
}
\date{November 28, 2026.}

\author{Shintaro~Akamine}
\address[Akamine]{%
College of Bioresource Sciences,
Nihon University, 1866 Kameino, Fujisawa, Kanagawa 252-0880, Japan
}
\email{akamine.shintaro@nihon-u.ac.jp}

\author{Masaaki Umehara}
\address[Umehara]{%
Department of Mathematical and Computing Sciences,
Institute of Science Tokyo,
2-12-1-W8-34, Ookayama, Meguro-ku,
Tokyo 152-8552, Japan
}
\email{umehara@c.titech.ac.jp}

\subjclass[2010]{Primary 53C42,\,\, Secondary 53C20,\,\, 53C50}

\thanks{%
The first author was partially supported by JSPS KAKENHI 
Grant Number 23K12979 and the
second author was partially supported by JSPS Grant-in-Aid
for Scientific Research (B), No. 23K20794.
}%

\begin{document}

\begin{abstract}
A flat Chebyshev metric defined on $\R^2$ is called an
``entire flat Chebyshev metric''. We prove that the plane with this metric
is isometric to the entire plane, an open half-plane,
or an open strip bounded by two parallel lines in Euclidean $2$-space.
Conversely, each of these domains, with its Euclidean metric, is isometric to $\R^2$
equipped with an entire flat Chebyshev metric.

As applications, we study time-like zero mean curvature (ZMC)
surfaces in Lorentz--Minkowski $3$-space. We characterize when
Magid's representation formula yields an entire time-like
ZMC graph over a space-like plane. We refer to such a graph
as a ``time-like ZMC graph of type $S$''.
We also give a new proof of a consequence of the Klotz--Milnor theorem,
which asserts that every member of the associated
one-parameter family of time-like ZMC graphs
of type $S$ is again an entire graph.
In contrast, we construct a time-like ZMC graph
over a space-like plane whose conjugate surface is not
an entire graph over any space-like plane.

We also use our completeness criterion to construct complete
flat surfaces in the unit sphere $S^3$ in $\R^4$
with unbounded mean curvature and entire
time-like ZMC graphs of type $S$ with unbounded Gaussian curvature.
\end{abstract}

\maketitle

\section*{Introduction}\label{sec0}
\setcounter{section}{0}

Flat Chebyshev metrics arise in the study of flat surfaces in the
$3$-sphere and time-like zero mean curvature surfaces in
Lorentz--Minkowski $3$-space.
We first define these metrics.

Throughout this paper, ``smooth'' means $C^\infty$-differentiable,
and $\R^2_{u,v}$ denotes the $uv$-plane. Let $\Omega_1$ and
$\Omega_2$ be smooth functions on $\R$ satisfying
\[
  \sin\bigl(\Omega_1(u)+\Omega_2(v)\bigr)>0
  \qquad (u,v\in\R).
\]
We call $(\Omega_1,\Omega_2)$ an \emph{$S$-pair}. It determines the
flat Riemannian metric
\begin{equation}\label{eq:1482}
  d\sigma^2
  :=du^2+2\cos\bigl(\Omega_1(u)+\Omega_2(v)\bigr)dudv+dv^2.
\end{equation}
We call $d\sigma^2$ a {\it flat Chebyshev metric}, which is said to 
be {\it entire} if it is defined on the whole of $\mathbb{R}^2$.

Entire flat Chebyshev metrics occur as the first fundamental
forms of the flat asymptotic Tchebychef immersions (FATs) into
the unit sphere $S^3$ in $\R^4$
considered by Kitagawa~\cite{Ki}.
Here we study the global geometry of these metrics without assuming
completeness. We will characterize the Euclidean plane domains that
admit global coordinates in which the metric is an entire flat
Chebyshev metric.
However, we first consider the case where the metric is complete:

\begin{definition}\label{def:1556}
An $S$-pair $(\Omega_1,\Omega_2)$
is \emph{complete} if the metric \eqref{eq:1482} is
complete. It is \emph{regular} if there exists $\epsilon>0$ such that
\begin{equation}\label{eq:1489}
  \sin\bigl(\Omega_1(u)+\Omega_2(v)\bigr)\geq\epsilon
  \qquad (u,v\in\R).
\end{equation}
We write $\mc P_S$ for the set of $S$-pairs, $\mc P_S^*$ for the
complete $S$-pairs, and $\mc P_S^{**}$ for the regular $S$-pairs.
\end{definition}

We have the following fact (cf. \cite{Ki,LMP}).

\begin{fact}\label{fact:227}
Every regular $S$-pair is complete, that is,
$\mc P_S^{**}\subset \mc P_S^*$.
\end{fact}

For a given $S$-pair $(\Omega_1,\Omega_2)$,
we define two
unit-speed plane curves by
\begin{equation}\label{eq:319}
  \gamma_1(u):=\int_0^u
  \begin{pmatrix}\cos\Omega_1(s)\\-\sin\Omega_1(s)\end{pmatrix}ds,
  \qquad
  \gamma_2(v):=\int_0^v
  \begin{pmatrix}\cos\Omega_2(s)\\\sin\Omega_2(s)\end{pmatrix}ds.
\end{equation}
We call them the \emph{developing curves} of the $S$-pair, and call
\begin{equation}\label{eq:329}
\Gamma:\R^2_{u,v}\ni (u,v) \mapsto \gamma_1(u)+\gamma_2(v)\in \R^2_{x,y}
\end{equation}
its \emph{developing map}.
This map is injective (cf. Proposition~\ref{prop:Inj}).
Since
\begin{equation}\label{eq:Jac}
  \det(\Gamma_u,\Gamma_v)
  =\sin\bigl(\Omega_1(u)+\Omega_2(v)\bigr)>0,
\end{equation}
and $\Gamma^*(dx^2+dy^2)=d\sigma^2$, the map $\Gamma$ is an
orientation-preserving local isometry.
Together with the injectivity of $\Gamma$, this gives the
following equivalences, which are related to Klotz--Milnor's
treatment of flat Chebyshev metrics~\cite{W}.

\begin{proposition}\label{prop:main}
For an $S$-pair $(\Omega_1,\Omega_2)$, the following are equivalent:
\begin{enumerate}
\item $(\Omega_1,\Omega_2)$ is complete;
\item $\Gamma:\R^2_{u,v}\to\R^2_{x,y}$ is a diffeomorphism;
\item $\Gamma:\R^2_{u,v}\to\R^2_{x,y}$ is proper;
\item $\Gamma(\R^2_{u,v})=\R^2_{x,y}$.
\end{enumerate}
\end{proposition}

\medskip
As an application, 
we obtain the following refinement of
Fact \ref{fact:227}.

\begin{corollary}\label{thm:K}
There exists an $S$-pair that is not complete.  Moreover, there
exists a complete $S$-pair that is not regular.  In particular,
$
  \mathcal P_S^{**}\subsetneq \mathcal P_S^*
  \subsetneq \mathcal P_S .
$
\end{corollary}

Viewing complete metrics as a special case of entire Chebyshev metrics, 
we are led to the following problem.

\medskip
\noindent
 {\it Characterize entire flat Chebyshev metrics
by the images of their developing maps.}

\medskip
The following theorem answers this question, including the
incomplete case.

\begin{introtheorem}\label{thm:A}
Let $D^2\subset\R^2$ be a nonempty simply connected domain. Then
$(D^2,dx^2+dy^2)$ is isometric to $(\R^2_{u,v},d\sigma^2)$ for some
entire flat Chebyshev metric $d\sigma^2$ if and only if $D^2$ is the
entire plane, an open half-plane, or an open strip bounded by two
parallel lines. Moreover, the developing map $\Gamma$
of every entire flat
Chebyshev metric maps $\R^2_{u,v}$
diffeomorphically onto a domain of one of these three types.
\end{introtheorem}

For the converse, we give explicit $S$-pairs for the plane, a
half-plane, and a strip of any prescribed width. For examples
of a strip and a half-plane, see \eqref{eq:IA1} and \eqref{eq:IA2},
respectively.

For an application of Theorem~\ref{thm:A}, we give the
following definition.

\begin{definition}\label{def:all}
A pair $(\gamma_1,\gamma_2)$ of simple regular plane curves is
\emph{unbounded in all directions} if, for every unit vector
$\mb e\in\R^2$, at least one of the functions
$\gamma_1\cdot\mb e$ and $\gamma_2\cdot\mb e$ is unbounded above.
Here ``$\cdot$'' denotes the Euclidean inner product.
\end{definition}

We then obtain the following practical criterion for the 
completeness of $S$-pairs.

\medskip
\noindent
{\bf Theorem B.}
{\it An $S$-pair is complete if and only if its pair of developing
curves is unbounded in all directions.}

\begin{remark*}
The $S$-pair assumption in Theorem~B is essential. The proper curves
$(s,s^2)$ and $(t,-t^2)$ are unbounded in all directions as a pair,
but their sum map omits $(0,1)$ since $s+t=0$ implies $s^2-t^2=0$.
\end{remark*}

We now turn to time-like zero mean curvature surfaces, abbreviated
to time-like ZMC surfaces. Let $\R^3_1$ be Lorentz--Minkowski
$3$-space with coordinates $(z,x,y)$ and metric
$-dz^2+dx^2+dy^2$. An immersion is time-like when its induced
metric is Lorentzian. By an \emph{entire graph over a plane}, we
mean a graph defined over the whole plane. We call graphs over
space-like, time-like, and light-like planes {\it types $S$, $T$, and
$L$}, respectively, and focus here on type $S$.

The following is the space-like plane case of Magid's
representation formula for time-like ZMC graphs.
For $\sigma\in\{1,-1\}$ and an $S$-pair $(\Omega_1,\Omega_2)$,
define a map $F:\R^2_{u,v}\to\R^3_1$, up to translation, by
\begin{equation}\label{eq:564}
  F_u(u,v)=(\sigma,\cos\Omega_1(u),-\sin\Omega_1(u)),
  \quad
  F_v(u,v)=(\sigma,\cos\Omega_2(v),\sin\Omega_2(v)).
\end{equation}
Magid \cite[Proposition~2']{M} showed that any time-like entire ZMC graph over the $xy$-plane is 
written as \eqref{eq:564}.
However,  if $F$ 
for a given $S$-pair $(\Omega_1,\Omega_2)$ 
gives an entire graph over the $xy$-plane,
it requires an additional condition, that is, 
the associated $S$-pair $(\Omega_1,\Omega_2)$ 
must be complete.
More precisely, we have the following proposition.

\begin{proposition}\label{cor:374}
The immersion $F$ as in \eqref{eq:564} gives an entire graph over the
$xy$-plane if and only if its $S$-pair is complete.
\end{proposition}

A direct computation shows that $F$ in \eqref{eq:564} has
zero mean curvature. Its Gaussian curvature is given by
(cf. Magid~\cite{M})
\begin{equation}\label{eq:K-Magid}
 K(u,v)=\frac{-\Omega_1'(u)\Omega_2'(v)}
 {\bigl(1-\cos(\Omega_1(u)+\Omega_2(v))\bigr)^2}
\qquad (u,v\in \R),
\end{equation}
where $\Omega'_1:=d\Omega_1/du$ and $\Omega'_2:=d\Omega_2/dv$.

Klotz--Milnor~\cite{W0,W} proved that
each member of the associated family
of a time-like entire harmonic graph is again an entire graph.
In terms of $S$-pairs, her theorem gives the following corollary.

\medskip
\noindent
{\bf Corollary C.}
{\it Let $(\Omega_1,\Omega_2)$ be a complete $S$-pair
$($resp.~an $S$-pair$)$.  Then
$$
  \bigl(\Omega_1(e^t u),\Omega_2(e^{-t}v)\bigr)
$$
is also a complete $S$-pair $($resp. an $S$-pair$)$
for each $t\in\R$.}

\medskip
Let $F$ be a time-like entire ZMC graph of type $S$.  Konderak's
representation formula~\cite{Kon} describes $F$ in terms of
Weierstrass data $(g,\omega)$.  For each $t\in\R$, the data
\[
  \bigl(g,(\cosh t+j\sinh t)\omega\bigr)
\]
produce a surface $F_t$.  Here the paracomplex numbers form the real
algebra
\[
  \check \C:=\{a+jb:a,b\in\R\},
  \qquad j^2=1,
\]
and $j$ is called the {\it paracomplex unit}.  The one-parameter family
$\{F_t\}_{t\in\R}$ is called the {\it associated family} of $F$.
If $(\Omega_1,\Omega_2)$ is the $S$-pair corresponding to $F$,
then the $S$-pair corresponding to $F_t$ is given by
$
  \bigl(\Omega_1(e^{-t}u),\Omega_2(e^t v)\bigr).
$
Thus Corollary~C implies that each member of the
associated family of a time-like entire ZMC graph of type $S$ is again
a time-like entire ZMC graph of type $S$.  In this paper, we give a new
proof of Corollary~C as an application of
Theorem~B.

A time-like ZMC surface $F$ also has a conjugate surface $F^*$,
obtained by taking the imaginary part of the representation formula.
The surface $F$ is obtained by taking the real part. For a space-like
ZMC surface, the conjugate surface belongs to its associated family.
In the time-like case, it does not. The Weierstrass data of $F^*$
are $(g,j\omega)$, and the first fundamental forms of $F$ and $F^*$
differ only by a sign (cf. \cite{Ac}).
As an application of Theorems~\ref{thm:A} and~B, we show that
the conjugate surface need not be an entire graph even when $F$ is
an entire graph.

\medskip
\noindent
{\bf Theorem D.}
{\it Let $F$ be an entire time-like ZMC graph over the $xy$-plane
$\Pi$, induced by a complete $S$-pair as in \eqref{eq:564}. Its
conjugate surface $F^*$ is a graph over a domain in $\Pi$ that is
the entire plane, an open half-plane, or an open strip
bounded by two parallel lines.

If the $S$-pair is regular, then $F^*$ is entire over $\Pi$.
However, there exists a complete non-regular $S$-pair for which
$F^*$ is not an entire graph over any space-like plane.}

\medskip
Krust's theorem states that if a minimal surface in Euclidean
$3$-space is represented as a graph over a convex domain in a plane,
then its conjugate surface is also a graph (see, for example,
\cite{Karcher}). An analogous result holds for space-like maximal
surfaces in $\R^3_1$: if such a surface is represented as a graph over
a convex domain in a space-like plane $\Pi$, then its conjugate surface
is also a graph over some domain in $\Pi$ (see \cite{AF,Lopez}). In
general, the domain of the conjugate graph need not coincide with the
original convex domain.
Theorem~D gives a time-like analogue of this result. Its final
assertion is stronger than saying that the conjugate surface is
not an entire graph over the original space-like plane.
In fact, in general,
a non-entire time-like ZMC graph over the $xy$-plane
can be an entire graph over another space-like plane
(cf. Example~\ref{ex:incomplete-over-another-plane}).

We next consider flat surfaces in $S^3$. 
The following fact is known.

\begin{fact}[cf. \cite{Ki}]\label{fact:B}
Let $(\Omega_1,\Omega_2)$ be a complete $S$-pair.  Then there
exists a complete flat immersion
$
\phi:\R^2_{u,v}\to S^3
$
whose first and second fundamental forms are given by
\eqref{eq:1482} and
\begin{equation}\label{eq:secondS3}
  2\sin\bigl(\Omega_1(u)+\Omega_2(v)\bigr)dudv,
\end{equation}
respectively.  Conversely, any complete flat immersion of
$\R^2_{u,v}$ into $S^3$ arises in this manner.
\end{fact}

Kitagawa~\cite[Theorem~4.3 and Remark~4.5]{Ki} showed that a complete
flat immersion $\phi:\R^2_{u,v}\to S^3$ associated with a
complete $S$-pair $(\Omega_1,\Omega_2)$ has a bounded mean curvature
function if and only if the $S$-pair is regular.
Consequently, Corollary~\ref{thm:K} yields the following.

\medskip
\noindent
{\bf Corollary E.}
{\it There exists a complete flat immersion $\phi:\R^2_{u,v}\to S^3$
whose mean curvature function is unbounded.}

\medskip
More precisely, the example
of a complete $S$-pair used in the proof of
Corollary~\ref{thm:K} induces a complete flat immersion
whose mean curvature function is unbounded both above and below
 (cf. Remark~\ref{rmk:917}).

Unlike the mean curvature of flat surfaces in $S^3$,
the Gaussian curvature of the associated time-like entire
ZMC graph depends on the first derivatives of the $S$-pair,
as in \eqref{eq:K-Magid}.
In particular, there exists a time-like entire ZMC graph of type $S$
with identically vanishing Gaussian curvature whose associated
$S$-pair is complete but not regular (see \eqref{eq:783}).
Thus, for time-like entire ZMC graphs, Gaussian curvature alone cannot
distinguish completeness from regularity of
the associated $S$-pair.
We prove the following theorem.

\medskip
\noindent
{\bf Theorem F.}
{\it There exists a time-like entire ZMC graph $F$ of type $S$ whose
associated $S$-pair $(\Omega_1,\Omega_2)$
is non-regular and has bounded derivatives
$\Omega_1'$ and $\Omega_2'$ on $\R$, whereas the Gaussian
curvature of $F$ is unbounded both above and below.}

\medskip
For fixed $u$ and $v$, respectively, the curves
\[
\R\ni u\mapsto \Gamma(u,v)\in\R^2_{x,y},
\qquad
\R\ni v\mapsto \Gamma(u,v)\in\R^2_{x,y}
\]
are integral curves of the two null direction fields induced by $F$.
Consequently, the boundedness of $\Omega_1'$ and $\Omega_2'$ is
equivalent to the uniform boundedness of the curvatures of these
integral curves in $\R^2_{x,y}$.
If the boundedness assumption on $\Omega_1'$ and $\Omega_2'$ is
dropped, there also exists a time-like entire ZMC graph of type $S$
whose Gaussian curvature is unbounded and whose associated $S$-pair
is regular  (see Example~\ref{ex:regular-unbounded-K}).
By contrast, the $S$-pair associated with the graph constructed in
Theorem~F is complete but not regular.

The existence of a complete flat embedding of the Euclidean plane
$(\mathbb{R}^2, dx^2+dy^2)$ into $S^3$ is a classical open 
problem (see \cite{AGM}).
Our numerical experiments indicate that the entire ZMC graph
used in the proof of Theorem~F does not give rise to an
embedded flat surface in $S^3$.

The paper is organized as follows.
In Section~\ref{sec:1}, we give basic properties of the
developing map and prove Proposition~\ref{prop:main} and
Corollary~\ref{thm:K}. Section~\ref{sec:2} gives properties
of $S$-pairs. In Section~\ref{sec:3}, we introduce
monotone coordinates and prove Theorems~\ref{thm:A} and~B.
Section~\ref{sec:4} gives applications to associated families,
conjugate surfaces, and Gaussian curvature.

\section{The developing map $\Gamma$ associated with an $S$-pair}\label{sec:1}

In this section, we study the developing map $\Gamma$ and relate
completeness of the flat Chebyshev metric to properties of its
developing curves.

\subsection*{Basic properties of the developing curves}

We fix an $S$-pair $(\Omega_1,\Omega_2)$ and consider the metric $d\sigma^2$
given by \eqref{eq:1482}.  Since
\[
  \sin(\Omega_1(u)+\Omega_2(v))>0 \qquad (u,v\in\R),
\]
and since $\R^2_{u,v}$ is connected, there exists an integer $k$ such that
\[
  2k\pi<\Omega_1(u)+\Omega_2(v)<(2k+1)\pi
  \qquad (u,v\in\R).
\]
Replacing $\Omega_2$ by $\Omega_2-2k\pi$ if necessary, we may and do assume
throughout this section that
\begin{equation}\label{eq:S-normal-section1}
  0<\Omega_1(u)+\Omega_2(v)<\pi
  \qquad (u,v\in\R).
\end{equation}
This replacement changes neither the metric $d\sigma^2$ nor the developing
map $\Gamma$.

We consider the smooth map
$
\Gamma:\R^2_{u,v}\to \R^2_{x,y}
$
defined by \eqref{eq:329}. 
For a matrix $M$, we denote its transpose by $M^T$.
If we write $\Gamma(u,v)=(x(u,v),y(u,v))^T$, then
\begin{equation}\label{eq:324b}
\Gamma_u(u,v)=
\begin{pmatrix}
\cos\Omega_1(u)\\
-\sin\Omega_1(u)
\end{pmatrix},
\qquad
\Gamma_v(u,v)=
\begin{pmatrix}
\cos\Omega_2(v)\\
\sin\Omega_2(v)
\end{pmatrix}.
\end{equation}
Hence
\begin{equation}\label{eq:Jac-section1}
\det(\Gamma_u,\Gamma_v)
=\sin(\Omega_1+\Omega_2)>0,
\end{equation}
so $\Gamma$ is an orientation-preserving local diffeomorphism.
Moreover,
\begin{equation}\label{eq:inner-products}
|\Gamma_u|^2=|\Gamma_v|^2=1,
\qquad
\Gamma_u\cdot \Gamma_v
=\cos(\Omega_1+\Omega_2),
\end{equation}
where $|\mb v|:=\sqrt{\mb v\cdot \mb v}$ for $\mb v\in \R^2_{x,y}$.
Therefore,
\begin{align}
\Gamma^*(dx^2+dy^2)
&=|\Gamma_u|^2du^2
+2(\Gamma_u\cdot \Gamma_v) dudv
+|\Gamma_v|^2dv^2 \notag\ \\
&=
du^2+2\cos(\Omega_1+\Omega_2)dudv+dv^2
=d\sigma^2.
\label{eq:pullback-flat}
\end{align}
We set
\[
\hat{\Omega}_1(u):=-\Omega_1(u)
\qquad (u\in\R).
\]
With this notation, the two curves can be written in the same form:
\[
\gamma_1(u):=
\int_0^u
\begin{pmatrix}
\cos\hat{\Omega}_1(\tilde u)\\
\sin\hat{\Omega}_1(\tilde u)
\end{pmatrix}
\,d\tilde u,
\qquad
\gamma_2(v):=
\int_0^v
\begin{pmatrix}
\cos\Omega_2(\tilde v)\\
\sin\Omega_2(\tilde v)
\end{pmatrix}
\,d\tilde v.
\]

\begin{proposition}\label{prop:bdd}
Let $(\Omega_1,\Omega_2)$ be an $S$-pair.
Then, for any $u,v\in \R$,
$\hat \Omega_1(u)\ne \Omega_2(v)$ holds. Moreover,
\begin{equation}\label{eq:O12bc}
|\hat \Omega_1(u_1)-\hat \Omega_1(u_2)|<\pi,\quad
|\Omega_2(v_1)-\Omega_2(v_2)|<\pi \qquad (u_1,u_2,v_1,v_2\in \R)
\end{equation}
also hold.
\end{proposition}

\begin{proof}
If $\hat\Omega_1(u)=\Omega_2(v)$, then
$\Omega_1(u)+\Omega_2(v)=0$, contradicting
\eqref{eq:S-normal-section1}.  This proves the first assertion.

For fixed $v$, both numbers
$\Omega_1(u_1)+\Omega_2(v)$ and
$\Omega_1(u_2)+\Omega_2(v)$ lie in $(0,\pi)$.  Their difference therefore
has absolute value strictly less than $\pi$, and hence
\[
|\Omega_1(u_1)-\Omega_1(u_2)|<\pi.
\]  The same argument, with $u$ fixed,
gives $|\Omega_2(v_1)-\Omega_2(v_2)|<\pi$.  Since
$\hat\Omega_1=-\Omega_1$, this is exactly \eqref{eq:O12bc}.
\end{proof}

The following statement follows immediately from \eqref{eq:O12bc}.

\begin{corollary}\label{cor:2}
Under the assumptions of Proposition~\ref{prop:bdd},
$\gamma_i'(s)$ does not coincide with $-\gamma_i'(t)$ for
$i=1,2$ and $s,t\in\mathbb{R}$ with $s\ne t$.
\end{corollary}

We next observe a similar relation between the two curves $\gamma_1$ and $\gamma_2$.

\begin{corollary}\label{cor:1}
Under the assumptions of Proposition~\ref{prop:bdd},
for any $(u,v)\in \mathbb{R}^2$,
$\gamma'_1(u)$ does not coincide with $\pm \gamma'_2(v)$.
\end{corollary}

\begin{proof}
Suppose first that
$\gamma'_1(u)=\gamma'_2(v)$.
Then $-\Omega_1(u)=\Omega_2(v)\pmod{2\pi}$ holds and hence
$\Omega_1(u)+\Omega_2(v)\in2\pi\mathbb Z$.
Suppose next that
$\gamma'_1(u)=-\gamma'_2(v)$.
Then $-\Omega_1(u)=\Omega_2(v)+\pi\pmod{2\pi}$ and hence
$\Omega_1(u)+\Omega_2(v)\in\pi+2\pi\mathbb Z$.  Both cases
contradict \eqref{eq:S-normal-section1}.
\end{proof}

We next show that the developing curves are simple and their sum map is injective.

\begin{proposition}\label{prop:439}
Each of the two curves $\gamma_i$ $(i=1,2)$ has no self-intersections.
\end{proposition}

\begin{proof}
If $\gamma_i$ contained a closed loop, then there would exist two points on the loop
such that the tangent vectors of $\gamma_i$ at these points are parallel 
but point in opposite directions,
contradicting Corollary~\ref{cor:2}.
\end{proof}

\begin{proposition}\label{prop:Inj}
Let $(\Omega_1,\Omega_2)$ be an $S$-pair.
Then its associated developing map $\Gamma$ is injective.
\end{proposition}

\begin{proof}
Suppose not. Then there exist two distinct points $(u_i,v_i)\in \mathbb{R}^2$ ($i=1,2$) such that
$\Gamma(u_1,v_1)=\Gamma(u_2,v_2)$, which implies
\[
\gamma_1(u_1)-\gamma_1(u_2)=\gamma_2(v_2)-\gamma_2(v_1).
\]
By Proposition~\ref{prop:439}, we must have $u_1\ne u_2$ and $v_1\ne v_2$.
Consider the line $\ell_1$ (resp.\ $\ell_2$) passing through the points
$\gamma_1(u_1),\gamma_1(u_2)$ (resp.\ $\gamma_2(v_1),\gamma_2(v_2)$).
Take a line parallel to $\ell_i$ and move 
it from infinity towards $\ell_i$ until it first meets 
the curve $\gamma_i$.
Since $\gamma_i$ has no self-intersections, 
the tangent to $\gamma_i$ at the point of 
intersection is parallel to $\ell_i$. 
Hence there exist $u_0$ between $u_1$ and $u_2$ 
and $v_0$ between $v_1$ and $v_2$ such that
$\gamma'_1(u_0)=\pm \gamma'_2(v_0)$, 
contradicting Corollary~\ref{cor:1}.
\end{proof}

We now prove Proposition~\ref{prop:main}.

\begin{proof}[Proof of Proposition~\ref{prop:main}]
$(1)\Rightarrow(2)$:
By Proposition~\ref{prop:Inj}, $\Gamma$ is injective.
By \eqref{eq:Jac}, $\Gamma$ is a local diffeomorphism.  Hence
$\Gamma$ is a diffeomorphism from $\R^2_{u,v}$ onto the open set
$\Gamma(\R^2_{u,v})\subset \R^2_{x,y}$.
Moreover, by \eqref{eq:pullback-flat},
$
  d\sigma^2=\Gamma^*(dx^2+dy^2).
$
Since $d\sigma^2$ is complete, the Euclidean metric on the open
set $\Gamma(\R^2_{u,v})$ is complete. A nonempty proper open subset
of the Euclidean plane is not complete: a straight segment from
an interior point to its first boundary point has finite length
and gives a Cauchy sequence with no limit in the domain.
Consequently, $\Gamma(\R^2_{u,v})=\R^2_{x,y}$.
Thus $\Gamma$ is a diffeomorphism from the $uv$-plane onto the
$xy$-plane.

$(2)\Rightarrow(1)$:  If $\Gamma$ is a diffeomorphism from the $uv$-plane onto
the $xy$-plane, then \eqref{eq:pullback-flat} shows that $d\sigma^2$ is the
pullback of the complete Euclidean metric $dx^2+dy^2$ on $\R^2_{x,y}$.  Therefore
$d\sigma^2$ is complete.

$(2)\Rightarrow(3)$:
Since $\Gamma$ is a homeomorphism from
$\R^2_{u,v}$ onto $\R^2_{x,y}$, it is proper.

$(3)\Rightarrow(4)$:  Suppose that $\Gamma$ is proper.  Since a proper map is
closed, $\Gamma(\R^2_{u,v})$ is a closed subset of $\R^2_{x,y}$.  Since $\Gamma$ is a
local diffeomorphism, $\Gamma(\R^2_{u,v})$ is also open in $\R^2_{x,y}$.  The image is
nonempty, and the target $\R^2_{x,y}$ is connected.  Hence
$\Gamma(\R^2_{u,v})=\R^2_{x,y}$.

$(4)\Rightarrow (2)$:  If $\Gamma(\R^2_{u,v})=\R^2_{x,y}$,
then Proposition~\ref{prop:Inj}
shows that $\Gamma$ is a bijective local diffeomorphism.  Hence $\Gamma$ is a
global diffeomorphism.
\end{proof}

\begin{proof}[Proof of Proposition~\ref{cor:374}]
Up to translation, \eqref{eq:564} gives
\[
 F(u,v)=\bigl(\sigma(u+v),\Gamma(u,v)\bigr).
\]
By Proposition~\ref{prop:Inj} and \eqref{eq:Jac}, $\Gamma$ is a
diffeomorphism onto its image. Hence $F$ is a graph over
$\Gamma(\R^2_{u,v})$, and it is entire over the $xy$-plane if and
only if $\Gamma(\R^2_{u,v})=\R^2_{x,y}$. By
Proposition~\ref{prop:main}, this is equivalent to completeness
of the $S$-pair.
\end{proof}

\begin{corollary}\label{cor:3}
Let $(\gamma_1,\gamma_2)$ be the pair of developing curves
associated with a complete $S$-pair $(\Omega_1,\Omega_2)$.
Then each $\gamma_i:\R\to\R^2_{x,y}$ $(i=1,2)$ is a proper map.
\end{corollary}

\begin{proof}
If $(\Omega_1,\Omega_2)$ is complete, then $\Gamma$ is proper by
Proposition~\ref{prop:main}.  Since
\[
  \gamma_1(t)=\Gamma(t,0),
  \qquad
  \gamma_2(t)=\Gamma(0,t),
\]
the inverse image $\gamma_1^{-1}(C)$ of a compact set $C\subset\R^2_{x,y}$
is identified with the intersection of the compact set $\Gamma^{-1}(C)$
with the line $v=0$. Therefore it is compact.  The argument for
$\gamma_2$ is identical.
\end{proof}

\subsection*{Basic properties of $S$-pairs}

The following lemma concerns Definition~\ref{def:all} and will
be used in Section~\ref{sec:3}.

\begin{lemma}\label{lem:linear-invariance}
Let $\gamma_i:\R\to\R^2_{x,y}$ $(i=1,2)$ be simple regular curves, and let
$\mc L:\R^2_{x,y}\to\R^2_{x,y}$ be an invertible linear map.  Then
$(\mc L\circ\gamma_1,\mc L\circ\gamma_2)$ is unbounded in all directions if and
only if $(\gamma_1,\gamma_2)$ is unbounded in all directions.
\end{lemma}

\begin{proof}
Let $\mb e$ be a unit vector in the target plane and let $\mc L^T$
denote the transpose of $\mc L$.  Since $\mc L$ is invertible,
$\mc L^T\mb e\ne \mb 0$, and
\[
 (\mc L\gamma_i(t))\cdot \mb e
 =\gamma_i(t)\cdot \mc L^T\mb e
 =|\mc L^T\mb e|\,
   \gamma_i(t)\cdot
   \frac{\mc L^T\mb e}{|\mc L^T\mb e|}.
\]
The scalar factor is positive, so the directional condition is preserved.
An invertible linear map and its inverse also preserve simplicity and regularity.  
Applying the same argument to $\mc L^{-1}$ proves
the converse.
\end{proof}

We use Proposition~\ref{prop:main} to prove Corollary~\ref{thm:K}.

\begin{proof}[Proof of Corollary~\ref{thm:K}]
We set
\[
\Omega_1(u):=\arcsin
\left(
\frac{1}{\sqrt{2\pi}}e^{-u^2/2}
\right),
\qquad
\Omega_2(v):=0
\qquad (u,v\in \R).
\]
Then $\Omega_1,\Omega_2$ are smooth functions on $\R$, and
$
        0<\Omega_1(u)+\Omega_2(v)<\pi
        \,\, (u,v\in\R).
$
Thus $(\Omega_1,\Omega_2)$ is an $S$-pair.
For this pair, the developing map is
\[
\Gamma(u,v)
=
\left(
v+\int_0^u \cos\Omega_1(t)\,dt,\,
-\int_0^u \sin\Omega_1(t)\,dt
\right)^T.
\]
Since
\[
\sin\Omega_1(t)=\frac{1}{\sqrt{2\pi}}e^{-t^2/2},
\]
the second component of $\Gamma$ takes values in
$
        \left(-\frac12,\frac12\right).
$
On the other hand, for each fixed $u$, the first component ranges over
the whole real line as $v$ varies.  Hence
\begin{equation}\label{eq:IA1}
        \Gamma(\R^2_{u,v})
        =
        \R\times\left(-\frac12,\frac12\right).
\end{equation}
Therefore, $\Gamma$ is not onto $\R^2_{x,y}$, and
the surface obtained from the representation formula is not a graph over
the entire $xy$-plane.  Equivalently, the metric $d\sigma^2$ is not
complete.  Thus $(\Omega_1,\Omega_2)\in\mc P_S\setminus\mc P_S^*$.

\medskip
We next consider the $S$-pair defined by
\begin{equation}\label{eq:783}
\Omega_1(u):=\frac{\pi}2+\arctan u,\qquad \Omega_2(v):=0.
\end{equation}
This pair is not regular.
We have
$$
\gamma_1(u)=\pmt{1-\sqrt{1+u^2}\\ -\op{arcsinh} u},\qquad
\gamma_2(v)=\pmt{v\\ 0},
$$
and
\[
\Gamma(u,v)=\pmt{v+1-\sqrt{1+u^2}\\ -\op{arcsinh}u}.
\]
This gives a diffeomorphism of the $uv$-plane onto the $xy$-plane.
In fact,
\[
\Phi:\R^2_{x,y}\ni \pmt{x\\ y} 
\mapsto \pmt{-\sinh y \\ x-1+\cosh y}\in \R^2_{u,v}
\]
is the inverse map of $\Gamma$. Thus $(\Omega_1,\Omega_2)$ is complete.
\end{proof}

\begin{remark} \label{rmk:917}
By Fact~\ref{fact:B}, the complete $S$-pair in
\eqref{eq:783} determines a complete flat immersion
$\phi:\R^2\to S^3$.  Put
$$
\Omega_1(u)+\Omega_2(v)
 =\frac{\pi}{2}+\arctan u.
$$
Since its first and second fundamental forms are
\eqref{eq:1482} and \eqref{eq:secondS3}, respectively, its
mean curvature is
\[
 H(u,v)=-\cot (\Omega_1(u)+\Omega_2(v)).
\]
Thus the flat surface associated with the complete $S$-pair in
\eqref{eq:783} has mean curvature $H(u,v)=u$, which is unbounded
both above and below.
\end{remark}

\section{Properties of $S$-pairs}\label{sec:2}

In this section, we study 
properties of $S$-pairs
without assuming completeness.
By
\eqref{eq:S-normal-section1},
we may assume that
\begin{equation}\label{eq:Sin2}
\epsilon \le \Omega_1(u)+\Omega_2(v) \le \pi-\epsilon \qquad (u,v\in\R),
\end{equation}
where $\epsilon\in [0,\pi/2)$.
If $\epsilon>0$, then $(\Omega_1,\Omega_2)$ is regular.
For $\epsilon=0$, no regularity condition is imposed.
Set
$\alpha:=\inf_{u\in\R}\Omega_1(u)$ and $\beta:=\sup_{u\in\R}\Omega_1(u)$.
Then we have $\alpha\le \Omega_1(0)\le\beta$.

\begin{proposition}\label{prop:588}
We have
\[
\epsilon-\alpha \le \Omega_2(v) \le \pi-\epsilon-\beta
\qquad (v\in\R).
\]
\end{proposition}

\begin{proof}
Fix $v_0\in\R$.
From \eqref{eq:Sin2}, we have
\[
\epsilon \le \Omega_1(u)+\Omega_2(v_0) \le \pi-\epsilon \qquad (u\in\R).
\]
The left inequality gives
$
\Omega_2(v_0) \ge \epsilon - \Omega_1(u) \,\, (u\in\R),
$
hence
\[
\Omega_2(v_0) \ge \sup_{u\in\R}(\epsilon-\Omega_1(u))
= \epsilon - \alpha.
\]
Similarly, the right inequality yields
$
\Omega_2(v_0) \le \pi-\epsilon-\Omega_1(u)\,\, (u\in\R),
$
hence
\[
\Omega_2(v_0) \le \inf_{u\in\R}\Big(\pi-\epsilon-\Omega_1(u)\Big)
= \pi-\epsilon-\beta.
\]
This holds for all $v_0\in\R$.
\end{proof}

By Proposition~\ref{prop:588},
we have
\begin{equation}\label{eq:627a}
\alpha \le \Omega_1 \le \beta, \qquad
\epsilon-\alpha \le \Omega_2 \le \pi-\epsilon-\beta.
\end{equation}

Two $S$-pairs $(\Omega_1,\Omega_2)$ and
$(\Omega_1^\star,\Omega_2^\star)$ are said to be \emph{$\R$-equivalent} if there
exists a constant $c\in\R$ such that
\begin{equation}\label{eq:720}
\Omega_1^\star := \Omega_1 + c,\qquad
\Omega_2^\star := \Omega_2 - c.
\end{equation}

\begin{proposition}\label{prop:714}
Let $(\Omega_1^\star,\Omega_2^\star)$ be $\R$-equivalent to an $S$-pair
$(\Omega_1,\Omega_2)$ satisfying \eqref{eq:Sin2}.
Then $(\Omega_1^\star,\Omega_2^\star)$
also satisfies \eqref{eq:Sin2}.
Moreover, $(\Omega_1^\star,\Omega_2^\star)$
is complete if and only if
$(\Omega_1,\Omega_2)$ is complete.
\end{proposition}

\begin{proof}
By \eqref{eq:720},
$\Omega_1^\star+\Omega_2^\star=\Omega_1+\Omega_2$.  Hence condition
\eqref{eq:Sin2} and the metric \eqref{eq:1482} are unchanged.  Moreover,
the tangent angles of both developing curves decrease by $c$, so their new
developing map is $R_{-c}\circ\Gamma$, where $R_{-c}$ denotes rotation by
angle $-c$ about the origin.  Thus completeness is preserved in both
directions.
\end{proof}

Using \eqref{eq:627a} and \eqref{eq:720}, we obtain
\[
-\beta-c \le -\Omega_1^\star \le -\alpha-c,\qquad
\epsilon-\alpha-c\le \Omega_2^\star\le \pi-\epsilon-\beta-c.
\]
We set $c:=-\beta-\epsilon$ and
$
  \mu:=\beta-\alpha+\epsilon
$.
By definition,
$
  \mu=-\inf \Omega_1^\star.
$
By \eqref{eq:627a},
\[
  0\le\mu=\beta-\alpha+\epsilon \le \pi-\epsilon\le\pi .
\]
In particular, $\mu\in[0,\pi]$.
Moreover, since
\[
\epsilon-\alpha-c=\epsilon-\alpha+(\beta+\epsilon)=\mu+\epsilon,
\]
we obtain
\begin{equation}\label{eq:733aaa}
  \epsilon \le -\Omega_1^\star \le \mu
  \le \mu+\epsilon \le \Omega_2^\star \le \pi .
\end{equation}
By setting $\epsilon=0$, the following statement is obtained immediately:

\begin{proposition}\label{prop:757}
Let $(\Omega_1,\Omega_2)$ be an $S$-pair satisfying  \eqref{eq:Sin2}.
Then there exist $\mu\in [0,\pi]$ and
an $S$-pair $(\Omega_1^\star,\Omega_2^\star)$
which is $\R$-equivalent to $(\Omega_1,\Omega_2)$
such that
\begin{equation}\label{eq:733ddd}
0 \le -\Omega_1^\star(u)\le \mu \le \Omega_2^\star(v)\le \pi.
\end{equation}
Moreover, if $(\Omega_1,\Omega_2)$ is complete, then
$(\Omega_1^\star,\Omega_2^\star)$ is also complete.
\end{proposition}

\section{Proof of Theorem~\ref{thm:A}}\label{sec:3}

Klotz--Milnor \cite[Theorem in Page 185]{W}
studied the curves $\gamma_1$ and $\gamma_2$ in \eqref{eq:319}
to prove a theorem underlying Corollary~C.
We also use these curves, and introduce monotone coordinates
to prove Theorem~\ref{thm:A}.
We first prove the following proposition.

\begin{proposition}\label{prop:821}
Let $(\gamma_1,\gamma_2)$ be the pair of developing curves
associated with an $S$-pair $(\Omega_1,\Omega_2)$.
If the $S$-pair $(\Omega_1,\Omega_2)$ is complete, then
$(\gamma_1,\gamma_2)$ is unbounded in all directions.
\end{proposition}

\begin{proof}
We fix a unit vector $\mb e$ in the $xy$-plane.
Suppose that $\gamma_1\cdot \mb e$ and $\gamma_2\cdot \mb e$ are both
bounded above.  Then $\Gamma\cdot\mb e$ is bounded above, which implies
$\Gamma(\R^2_{u,v})\ne\R^2_{x,y}$.  This contradicts completeness by
Proposition~\ref{prop:main}.  Hence the pair is unbounded in all directions in the
sense of Definition~\ref{def:all}.
\end{proof}

For the moment, we assume that neither $\gamma_1$ nor
$\gamma_2$ is a straight line.
Since $\sin(\Omega_1+\Omega_2)>0$ on the connected plane, there is an
integer $k$ such that
$$
2k\pi<\Omega_1+\Omega_2<(2k+1)\pi.
$$  
Replacing $\Omega_2$ by
$\Omega_2-2k\pi$ preserves the metric $d\sigma^2$ in \eqref{eq:1482}
and the developing curves in \eqref{eq:319}, and gives
$0<\Omega_1+\Omega_2<\pi$.  Thus
\eqref{eq:Sin2} holds with $\epsilon=0$, so Proposition~\ref{prop:757}
applies.
We may therefore further assume that
the given $S$-pair satisfies
\[
0 \le \hat \Omega_1 \le  \mu \le \Omega_2 \le \pi.
\]
If $\mu=\pi$ (resp. $\mu=0$),
then $\Omega_2$ (resp. $\Omega_1$) is constant, and hence
$\gamma_2$ (resp. $\gamma_1$)
is a straight line, contradicting our assumption.
Therefore, we may assume that $0<\mu<\pi$.
Thus, we work under the following assumption:
\begin{equation}\tag{$\diamond$}
0 \le \hat\Omega_1 \le \mu \le \Omega_2 \le \pi
\qquad (0<\mu<\pi).
\end{equation}

\begin{lemma}\label{lem:developing-curves-proper}
Assume $(\diamond)$.  Then both developing curves
$\gamma_1,\gamma_2:\R\to\R^2_{x,y}$ are proper.
\end{lemma}

\begin{proof}
For $\theta\in\R$, set
$\mb e_\theta=(\cos\theta,\sin\theta)$.
Since the tangent angle of $\gamma_1$ is $\hat\Omega_1$, we have
\[
 \frac{d}{du}
 \bigl(\gamma_1(u)\cdot\mb e_{\mu/2}\bigr)
 =\cos\left(\hat\Omega_1(u)-\frac{\mu}{2}\right)
 \ge \cos\frac{\mu}{2}>0.
\]
It follows that
$\gamma_1(u)\cdot\mb e_{\mu/2}\to \pm \infty$ as
$u\to\pm \infty$.
Hence $\gamma_1$ is proper.
Similarly,
\[
 \frac{d}{dv}
 \bigl(\gamma_2(v)\cdot\mb e_{(\mu+\pi)/2}\bigr)
 =\cos\left(\Omega_2(v)-\frac{\mu+\pi}{2}\right)
 \ge \cos\frac{\pi-\mu}{2}>0.
\]
Hence $\gamma_2(v)\cdot\mb e_{(\mu+\pi)/2} \to \pm \infty$ as
$v\to\pm \infty$.
Hence $\gamma_2$ is also proper.
\end{proof}

Without $(\diamond)$, $\gamma_1$ and $\gamma_2$ need not both be proper:

\begin{example}\label{ex:straight-line-nonproper}
Consider the simple regular curve
$
  c(t):=(\sin t,-\arctan t)^T
  \,\,(t\in\R).
$
Its speed is
\[
  |c'(t)|
  =
  \sqrt{\cos^2t+\frac{1}{(1+t^2)^2}}>0
\]
and the length of $c$ is infinite in both directions. Hence
we can write
$$
  \gamma_1(u):=c(t(u))
  \qquad (u\in\R),
$$
where $u$ is an arc-length parameter of $c$. Then
$\dot\gamma_1:=d\gamma_1/du$ satisfies
\[
\dot\gamma_1(u)
  =
  \frac{1}{
  \sqrt{\cos^2t+(1+t^2)^{-2}}
  }
  \pmt{
    \cos t\\ -\frac{1}{1+t^2}
  }
  \qquad t=t(u).
\]

Let
$
  \operatorname{arccot}:\R\rightarrow(0,\pi)
$
denote the inverse of the cotangent function restricted to
$(0,\pi)$, and define
\[
  \Omega_1(u):=
  \operatorname{arccot}
  \Big(
    (1+t(u)^2)\cos t(u)
  \Big),
  \qquad
  \Omega_2(v):=0.
\]
In particular, $(\diamond)$ fails for this $S$-pair.
In this setting,  we have
\[
  \cos\Omega_1(u)
  =
  \frac{(1+t^2)\cos t}
  {\sqrt{1+(1+t^2)^2\cos^2t}},\qquad
  \sin\Omega_1(u)
  =
  \frac{1}
  {\sqrt{1+(1+t^2)^2\cos^2t}}.
\]
Consequently,
\[
\dot\gamma_1(u)
  =
  \pmt{
    \cos\Omega_1(u)\\-\sin\Omega_1(u)
  }.
\]
Since
$
  0<\Omega_1(u)<\pi,
$
we have
\[
  \sin\bigl(\Omega_1(u)+\Omega_2(v)\bigr)
  =
  \sin\Omega_1(u)>0
  \qquad ((u,v)\in\R^2).
\]
Thus $(\Omega_1,\Omega_2)$ is an $S$-pair.
The corresponding second developing curve is
\[
  \gamma_2(v)
  =
  \int_0^v
  \begin{pmatrix}
    \cos\Omega_2(\tilde v)\\
    \sin\Omega_2(\tilde v)
  \end{pmatrix}
  d\tilde v
  =\pmt{v\\ 0},
\]
and hence $\gamma_2$ is a straight line.
Since
\[
  \gamma_1(\R)
  =
  c(\R)
  \subset A:=
  [-1,1]\times
  \left[-\frac{\pi}{2},\frac{\pi}{2}\right],
\]
we have
$
  \gamma_1^{-1}(A)=\R,
$
which is not compact. Hence $\gamma_1$ is not proper.

\end{example}

\begin{proposition}\label{prop:934b}
Let $\gamma:\R\to\R^2$ be a simple regular curve that is proper.
Then $\gamma$ separates $\R^2$ into two domains.
\end{proposition}

\begin{proof}
The one-point compactification of $\R$ is $S^1$, and that of
$\R^2$ is $S^2$.
Since $\gamma$ is proper, it extends continuously to a map
$\tilde\gamma:S^1\to S^2$.
Since $\gamma$ is simple, $\tilde\gamma$ is a simple closed curve in $S^2$.
The Jordan curve theorem gives two complementary domains.
In the one-point compactification, the point at infinity belongs to the
Jordan curve and hence to the boundary of both complementary components.
Therefore both components are unbounded in $\R^2$.
\end{proof}

\subsection*{Monotone coordinate systems}

In this subsection, we introduce the monotone coordinate system
associated with an $S$-pair $(\Omega_1,\Omega_2)$ under 
the condition $(\diamond)$.
Set
\[
\mathbf e_1:=\pmt{1\\ 0},\qquad
 \mathbf e_2:=\pmt{\cos\mu\\ \sin\mu}.
\]
Since $0<\mu<\pi$, we have
$
 \det(\mathbf e_1,\mathbf e_2)=\sin\mu>0.
$
Thus $(\mathbf e_1,\mathbf e_2)$ is a basis of
$\R^2_{x,y}$.
Let
$
 \mathcal L_\mu:\R^2_{x,y}\rightarrow \R^2_{\xi,\eta}
$
be the linear isomorphism such that
\[
 \mathcal L_\mu(\xi\mathbf e_1+\eta\mathbf e_2)
 =\pmt{\xi\\ \eta}.
\]
We set
\begin{equation}\label{eq:1234}
 \bar\gamma_i
 :=\mathcal L_\mu\circ\gamma_i
 =\pmt{\xi_i\\ \eta_i}
 \quad (i=1,2),\qquad
 \bar\Gamma
 :=\mathcal L_\mu\circ\Gamma
 =\bar\gamma_1+\bar\gamma_2.
\end{equation}

\begin{lemma}\label{prop:701}
For every $u,v\in\R$,
\[
 \xi'_1(u)\ge0,\qquad \eta'_1(u)\ge0,
 \qquad
 \xi'_2(v)\le0,\qquad \eta'_2(v)\ge0.
\]
\end{lemma}

\begin{proof}

By $(\diamond)$, we have $0 \le \hat\Omega_1 \le \mu$
and $\mu \le \Omega_2 \le \pi$, where $0<\mu<\pi$.
Hence the vectors
\[
 \gamma'_1(u)=\pmt{\cos\hat\Omega_1(u)\\ \sin\hat\Omega_1(u)},\qquad
 \gamma'_2(v)=\pmt{\cos\Omega_2(v)\\ \sin\Omega_2(v)}
\]
lie in the closed sectors spanned by $\mathbf e_1,\mathbf e_2$ and
by $\mathbf e_2,-\mathbf e_1$, respectively.
Applying
$\mathcal L_\mu$ gives the required signs.
\end{proof}

\subsection*{Monotone pairs of curves in the $\xi\eta$-plane}

We define monotone pairs of curves in the $\xi\eta$-plane and
study their sum map $G$. We determine its image under suitable
unboundedness assumptions.

\begin{definition}
Let
$c_i\colon \mathbb{R}\rightarrow \mathbb{R}^2_{\xi,\eta}$
$(i=1,2)$
be smooth regular curves. The pair $(c_1,c_2)$ is called a
\emph{monotone pair} if the following conditions are satisfied:
\begin{enumerate}
\item $c_1(0)=c_2(0)$,
\item each map $c_i=(\xi_i,\eta_i)^T$ has no self-intersections,
\item each map $c_i=(\xi_i,\eta_i)^T$ is a proper map 
into $\mathbb{R}^2_{\xi,\eta}$,
\item the velocity vector
$
c_1'=(\xi_1',\eta_1')^T
$
lies in the closure of the first quadrant of
$\mathbb{R}^2_{\xi,\eta}$, that is,
$\xi_1',\eta_1'\geq 0$, and
\item the velocity vector
$
c_2'=(\xi_2',\eta_2')^T
$
lies in the closure of the second quadrant of
$\mathbb{R}^2_{\xi,\eta}$, that is,
$\xi_2'\leq 0$ and
$\eta_2'\geq 0$.
\end{enumerate}
For such a pair $(c_1,c_2)$, we set
\[
G\colon \mathbb{R}^2_{u,v}\ni (u,v)
\mapsto c_1(u)+c_2(v)\in \mathbb{R}^2_{\xi,\eta}.
\]
We call $G$ the \emph{developing map associated with} $(c_1,c_2)$.
Straight lines are allowed in this definition.
\end{definition}

The following proposition follows from
Lemmas~\ref{lem:developing-curves-proper} and~\ref{prop:701}.

\begin{proposition}\label{prop:1316}
Let $\gamma_1$ and $\gamma_2$ be the \emph{developing curves}
associated with an $S$-pair $(\Omega_1,\Omega_2)$ satisfying
$(\diamond)$.
Then $(\bar \gamma_1,\bar \gamma_2)$ and
$\bar \Gamma$ given in \eqref{eq:1234}
give a monotone pair and its developing map, respectively.
\end{proposition}

\begin{definition}\label{def:4dir}
Let $(c_1,c_2)$ be a monotone pair, where
$c_i=(\xi_i,\eta_i)^T$.
Fix $i\in \{1,2\}$.
We say that $c_i$ is
\emph{unbounded to the left}, \emph{unbounded to the right},
\emph{unbounded below}, or \emph{unbounded above} if
$$
\begin{array}{@{}l@{}}
\inf(\xi_i(\R))=-\infty,\,\,\,
\sup(\xi_i(\R))=+\infty,\,\,\,
\inf(\eta_i(\R))=-\infty,\,\, \text{or} \,\,
\sup(\eta_i(\R))=+\infty,
\end{array}
$$
respectively.
Moreover,
\begin{itemize}
\item
$(c_1,c_2)$ is
\emph{unbounded to the left}, \emph{unbounded to the right},
\emph{unbounded below}, or \emph{unbounded above} if
at least one of $c_1$ and $c_2$ has the corresponding property,
\item
$(c_1,c_2)$ is said to be
\emph{bounded from below} if it is not
unbounded below,
\item
$(c_1,c_2)$ is
\emph{unbounded in the four coordinate directions} if it is
unbounded to the left, unbounded to the right, unbounded below, and
unbounded above.
\end{itemize}
\end{definition}

\begin{proposition}\label{prop:all-to-four}
If a monotone pair $(c_1,c_2)$ is unbounded in all directions in $\R^2_{\xi,\eta}$,
then $(c_1,c_2)$ is unbounded in the four
coordinate directions in $\R^2_{\xi,\eta}$.
\end{proposition}

\begin{proof}
Apply Definition~\ref{def:all} to the four unit vectors
$$
\pmt{1\\0}, \quad \pmt{-1\\ 0},\quad \pmt{0\\ 1},\quad \pmt{0\\ -1}.
$$  The resulting four
conditions are precisely those in Definition~\ref{def:4dir}.
\end{proof}

Let $\mc R$ denote the linear transformation of the $\xi\eta$-plane
given by the counterclockwise rotation through $90^\circ$ about the origin.
Set $\iota(t):=-t$ for $t\in\R$.
Then the following statement also follows immediately from the definition.

\begin{proposition}\label{prop:1331}
If $(c_1,c_2)$ is a monotone pair,
then
\[
\mc R(c_1,c_2):=(\mc R\circ c_2\circ\iota,\mc R \circ c_1)
\] is also a monotone pair.
Moreover, if $G$ and $\tilde G$ are the developing maps of
$(c_1,c_2)$ and $\mc R(c_1,c_2)$, respectively, then
$
\tilde G(u,v)=\mc R\circ G(v,-u).
$
\end{proposition}

\begin{proof}
Put $\tilde c_1:=\mc R\circ c_2\circ\iota$ and
$\tilde c_2:=\mc R\circ c_1$.  Since
$\mc R(x,y)^T=(-y,x)^T$, we have
\[
 \tilde c_1'(u)
 =\pmt{\eta_2'(-u)\\-\xi_2'(-u)},
 \qquad
 \tilde c_2'(v)
 =\pmt{-\eta_1'(v)\\ \xi_1'(v)}.
\]
Thus $\tilde c_1'$ lies in the closure of the first quadrant and
$\tilde c_2'$ lies in the closure of the second quadrant.
Moreover, $\tilde c_1(0)=\tilde c_2(0)$, and regularity and
properness are preserved.  Hence $(\tilde c_1,\tilde c_2)$ is a
monotone pair.  The formula for $\tilde G$ follows directly from the
definition.
\end{proof}

As a consequence, we obtain the following.

\begin{corollary}\label{cor:1369}
Let $(c_1,c_2)$ be a monotone pair.
If it is not unbounded in the four
coordinate directions, then there exists $i\in \{0,1,2,3\}$
such that $\mc R^i(c_1,c_2)$ is bounded from below,
where $\mc R^i$ denotes the $i$-fold iteration of the operation in
Proposition~\ref{prop:1331}.
\end{corollary}

\begin{proof}
At least one of the four conditions in Definition~\ref{def:4dir}
fails.  Choose $i\in\{0,1,2,3\}$ so that rotation through $90i$
degrees sends the corresponding missing direction to the downward
coordinate direction.  The interchange and reparametrization in
Proposition~\ref{prop:1331} do not affect boundedness in a fixed
direction.  Thus, for the transformed pair, neither curve is
unbounded downward, that is, both of its $\eta$-coordinates are
bounded from below.
\end{proof}

\subsection*{Proofs of Theorems~\ref{thm:A} and~B}

In connection with Corollary~\ref{cor:1369}, we prove the following.

\begin{proposition}\label{prop:1385}
Let $(c_1,c_2)$ be a monotone pair, and let $G$ be its developing map.
Suppose that $(c_1,c_2)$ is bounded from below.
Then $c_1$ is unbounded to the left and
$c_2$ is unbounded to the right.
\end{proposition}

\begin{proof}
Write $c_i=(\xi_i,\eta_i)^T$.  For $u\leq0$, monotonicity gives
$$
 \xi_1(u)\leq\xi_1(0),\qquad \eta_1(u)\leq\eta_1(0).
$$
If $\xi_1$ were bounded from below, then the lower boundedness of
$\eta_1$ would imply that $c_1((-\infty,0])$ is contained in a
compact rectangle, contradicting the properness of $c_1$.  Hence
$\inf(\xi_1(\R))=-\infty$.
Similarly, for $v\leq0$,
$$
 \xi_2(v)\geq\xi_2(0),\qquad \eta_2(v)\leq\eta_2(0).
$$
If $\xi_2$ were bounded from above, then $c_2((-\infty,0])$ would be
contained in a compact rectangle, contradicting the properness of
$c_2$.  Thus $\sup(\xi_2(\R))=+\infty$.
\end{proof}

Under these assumptions, the image contains the following horizontal lines.

\begin{proposition}\label{prop:1403b}
Let $(c_1,c_2)$ be a monotone pair, and let $G$ be its developing map.
Suppose that
$c_1$ is unbounded to the left and
$c_2$ is unbounded to the right.
If $(a,b)^T=G(u_0,v_0)$,  then the horizontal line
\[
L_b:=\left\{\pmt{\xi\\ b}\in  \R^2_{\xi,\eta}\,;\, \xi\in \R\right\}
\]
is a subset of $G(\R^2_{u,v})$.
\end{proposition}

We first prove the following lemma.

\begin{lemma}\label{lem:1403b}

Let $(c_1,c_2)$ be a monotone pair such that $c_1$ is unbounded to the
left, $c_2$ is unbounded to the right, and
$c_1(0)=c_2(0)=\pt{O}$, where $\pt{O}$ is the origin of the $\xi\eta$-plane.  
Then the horizontal line
\[
L_0:=\left\{\pmt{\xi\\ 0}\in\R^2_{\xi,\eta}\,;\,\xi\in\R\right\}
\]
is a subset of $G(\R^2_{u,v})$.

\end{lemma}

\begin{figure}[ht!]
  \begin{center}
    \includegraphics[width=.48\linewidth,height=56mm,keepaspectratio]{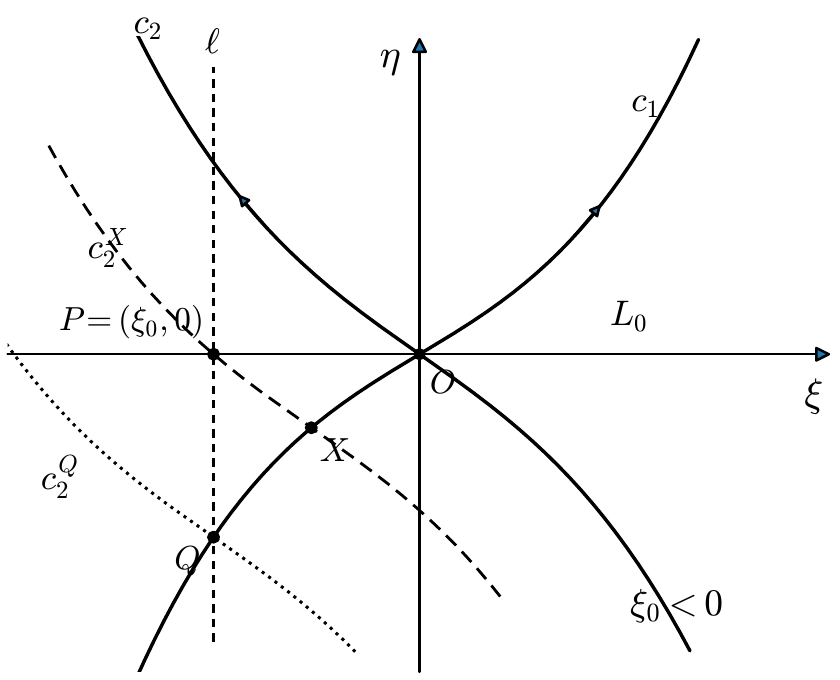}\,\,
    \includegraphics[width=.48\linewidth,height=56mm,keepaspectratio]{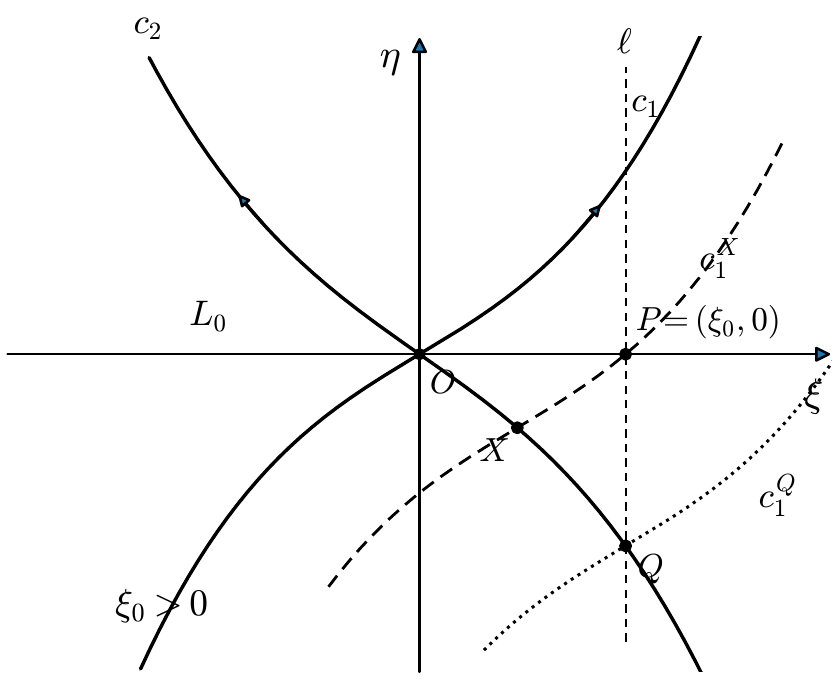}
    \caption{
The case $\xi_0<0$ (left) and the case $\xi_0>0$ (right).
}
    \label{fig12}
  \end{center}
\end{figure}

\begin{proof}
Fix $\pt{P}:=(\xi_0,0)\in L_0$.  If $\xi_0=0$, then $\pt{P}=G(0,0)$.
Suppose first that $\xi_0<0$, and let $\ell$ be the line parallel to the
$\eta$-axis through $\pt{P}$.  Since $c_1$ is unbounded to the left, $\ell$
meets $c_1$ at a point $\pt{Q}=c_1(a)$ for some $a<0$.  For
$\pt{X}\in c_1([a,0])$, set
\[
 c_2^{\pt{X}}(v):=\pt{X}+c_2(v).
\]
We now use
Proposition \ref{prop:934b}.
By monotonicity, $\pt{P}$ lies in the closure of the lower domain bounded by
$c_2^{\pt{O}}=c_2$ and in the closure of the upper 
domain bounded by $c_2^\pt{Q}$.
By continuity of this translated family, there exist $u_0\in[a,0]$ and
$v_1\in\R$ such that (cf. Figure~\ref{fig12}, left)
\[
\pt{P}=c_2^{c_1(u_0)}(v_1)=c_1(u_0)+c_2(v_1)=G(u_0,v_1).
\] 

Suppose next that $\xi_0>0$.  Since $c_2$ is unbounded to the right,
$\ell$ meets $c_2$ at a point $\pt{Q}=c_2(b)$ for some $b<0$.  For
$\pt{X}\in c_2([b,0])$, set $c_1^\pt{X}(u):=\pt{X}+c_1(u)$.
The same argument gives $v_0\in[b,0]$ and $u_1\in\R$ such that 
(cf. Figure~\ref{fig12}, right)
\[
 \pt{P}=c_1^{c_2(v_0)}(u_1)=c_1(u_1)+c_2(v_0)=G(u_1,v_0).
\]
Thus $L_0\subset G(\R^2_{u,v})$.
\end{proof}

\begin{proof}[Proof of Proposition~\ref{prop:1403b}]

Define
\[
 \check c_1(u):=c_1(u+u_0)-c_1(u_0),\qquad
 \check c_2(v):=c_2(v+v_0)-c_2(v_0).
\]
Then $(\check c_1,\check c_2)$ is a monotone pair satisfying
$\check c_1(0)=\check c_2(0)=\pt{O}$, and it has the same left- and
right-unboundedness properties as $(c_1,c_2)$.  Its developing map is
\[
 \check G(u,v)
 =G(u+u_0,v+v_0)-\pmt{a\\ b}.
\]
By Lemma~\ref{lem:1403b}, $L_0\subset\check G(\R^2)$.  Consequently,
\[
 L_b=\pmt{a\\ b}+L_0\subset G(\R^2),
\]
which proves the assertion.
\end{proof}

We now prove the following theorem.

\begin{theorem}\label{prop:1386}
Let $(c_1,c_2)$ be a monotone pair, and let $G$ be its developing map.
Suppose that
$c_1$ is unbounded to the left and
$c_2$ is unbounded to the right.
Assume, in addition, that $G$ is a local diffeomorphism.
Then
$
G(\R^2_{u,v})=\R\times (\alpha,\beta),
$
where
\begin{align*}
&\alpha:=\inf_{u,v\in \R}(\eta_1(u)+\eta_2(v))\in \R\cup\{-\infty\}, \\
&\beta:=\sup_{u,v\in \R}(\eta_1(u)+\eta_2(v))\in \R\cup\{+\infty\}.
\end{align*}
\end{theorem}

\begin{proof}
We set $I := (\alpha,\beta)$ and
\[
 h(u,v):=\eta_1(u)+\eta_2(v).
\]

Since $G$ is a local diffeomorphism, $dh$ never vanishes: otherwise
$\eta_1'(u)=\eta_2'(v)=0$, which would imply
$\det(c_1'(u),c_2'(v))=0$.  Hence $h$ is a submersion, and
$h(\R^2)$ is open.  Since $\R^2$ is connected, $h(\R^2)$ is an open
interval. By the definitions of $\alpha$ and $\beta$, it is precisely
$I$.
Consequently,
\begin{equation}\label{eq:1514}
G(\mathbb{R}_{u,v}^2)\subset \mathbb{R}\times I.
\end{equation}
It therefore suffices to prove the reverse inclusion.

Take an arbitrary $b\in I$. Then there exist $a\in\mathbb{R}$ and
$(u_0,v_0)\in\mathbb{R}_{u,v}^2$ such that
$
(a,b)^T=G(u_0,v_0).
$
By Proposition~\ref{prop:1403b}, the horizontal line
\[
L_b:=\left\{\pmt{\xi\\ b}\in\mathbb{R}_{\xi,\eta}^2 \,;\, \xi\in\mathbb{R}\right\}
\]
is contained in $G(\mathbb{R}_{u,v}^2)$.
Since $b\in I$ was arbitrary, we obtain
\[
\mathbb{R}\times I
=\bigcup_{b\in I}L_b
\subset G(\mathbb{R}_{u,v}^2).
\]
Combining this with \eqref{eq:1514}, we conclude that
$G(\mathbb{R}_{u,v}^2)=\mathbb{R}\times I$.
This completes the proof.
\end{proof}

\begin{corollary}\label{cor:1436}
Let $\gamma_1$ and $\gamma_2$ be the \emph{developing curves}
associated with an $S$-pair $(\Omega_1,\Omega_2)$ satisfying
$(\diamond)$.
If $(\bar \gamma_1,\bar \gamma_2)$ in \eqref{eq:1234} is
unbounded in the four
coordinate directions,
then $(\Omega_1,\Omega_2)$ is a complete $S$-pair.
\end{corollary}

\begin{proof}
Write $c_i=(\xi_i,\eta_i)^T:=\bar\gamma_i$ ($i=1,2$)
and
$G:=\bar\Gamma=c_1+c_2$. We first prove that $G$ is proper.
Otherwise, there would exist a sequence $(u_n,v_n)$ with
$|u_n|+|v_n|\to\infty$ such that $G(u_n,v_n)$ stays in a compact
set. Since both $c_1$ and $c_2$ are proper, after passing to a
subsequence each of $u_n$ and $v_n$ tends to $+\infty$ or
$-\infty$. Indeed, if one parameter stayed bounded while the
other diverged, the corresponding sum could not stay bounded.

Suppose first that $u_n,v_n\to+\infty$. Both
$\eta_1(u_n)$ and $\eta_2(v_n)$ are bounded below, and their sum
is bounded, so each sequence is bounded above. Monotonicity then
implies that both $\eta_1$ and $\eta_2$ are bounded above on
$\R$. This contradicts unboundedness above. If
$u_n,v_n\to-\infty$, the same argument gives a common lower
bound for $\eta_1$ and $\eta_2$, contradicting unboundedness
below.

If $u_n\to+\infty$ and $v_n\to-\infty$, both
$\xi_1(u_n)$ and $\xi_2(v_n)$ are bounded below. Their bounded
sum forces both $\xi_1$ and $\xi_2$ to be bounded above on
$\R$, contradicting unboundedness to the right. Finally, if
$u_n\to-\infty$ and $v_n\to+\infty$, both $\xi_1$ and
$\xi_2$ are bounded below on $\R$, contradicting unboundedness
to the left. These four cases exhaust all possibilities.

Thus $G$ is proper. Since $\mathcal L_\mu$ is an invertible
linear map, $\Gamma$ is also proper. Proposition~\ref{prop:main}
now proves completeness.
\end{proof}

\begin{proof}[Proof of Theorem~\ref{thm:A}]
Let us assume that $(D^2,dx^2+dy^2)$ is isometric to 
$(\R^2_{u,v},d\sigma^2)$ for some
entire flat Chebyshev metric $d\sigma^2$ in \eqref{eq:1482}. 
Then  $(\Omega_1,\Omega_2)$ gives an $S$-pair.
We first prove that $\Gamma$ in \eqref{eq:329} gives an isometry onto one of 
three domains in the assertion.
The injectivity of $\Gamma$ follows from Proposition~\ref{prop:Inj}.
We show the assertion concerning its image.

If $\gamma_1$ is a straight line, then, after rotating the target
plane, we may assume that
\[
  \gamma_1(u)=\pmt{u\\ 0}.
\]
Indeed, $\gamma_1$ has unit speed and $\gamma_1(0)=\pt{O}$, and its unit
tangent vector is constant.  Write
$
  \gamma_2(v)=(x_2(v),y_2(v))^T.
$
Since $\Gamma$ is a local diffeomorphism, we have
\[
  0\ne \det(\Gamma_u,\Gamma_v)=y_2'(v).
\]
Consequently, $y_2$ is a local diffeomorphism, and
$
  I:=y_2(\R)
$
is a nonempty open interval.  Moreover,
\[
  \Gamma(\R^2_{u,v})
  =\left\{\pmt{u+x_2(v)\\ y_2(v)}\,;\ u,v\in\R\right\}
  =\R\times I.
\]
According as $I$ is bounded, unbounded at exactly one end, or equal to
$\R$, this image is an open strip, an open half-plane, or the
entire plane.  The case where $\gamma_2$ is a straight line is proved
in the same way.  Thus we may assume that neither $\gamma_1$ nor
$\gamma_2$ is a straight line.
We may therefore assume that $(\Omega_1,\Omega_2)$ satisfies $(\diamond)$.  By
Proposition~\ref{prop:1316}, a monotone pair
$(\bar\gamma_1,\bar\gamma_2)$ in $\R^2_{\xi,\eta}$ is induced.

If $(\bar\gamma_1,\bar\gamma_2)$ is unbounded in the four coordinate
directions, Corollary~\ref{cor:1436} and Proposition~\ref{prop:main}
 show that
$\Gamma(\R^2_{u,v})=\R^2_{x,y}$.
Otherwise, Corollary~\ref{cor:1369} gives an
$i\in\{0,1,2,3\}$ such that
\[
(c_1, c_2):=\mc R^i(\bar\gamma_1,\bar\gamma_2)
\]
is bounded from below. By Proposition~\ref{prop:1331}, the developing map
 associated with $(c_1,c_2)$ has image
$$
\mc R^i\circ \bar\Gamma(\R^2_{u,v})
$$
and is a local diffeomorphism. Proposition~\ref{prop:1385} shows that $c_1$ 
is unbounded to the left and that $c_2$ is unbounded to the right. 
Hence, by Theorem~\ref{prop:1386},
$
\mc R^i\circ \bar\Gamma(\R^2_{u,v})
$
is either an open horizontal strip or an open half-plane.
Applying the inverse rotation and then $\mathcal L_\mu^{-1}$,
we conclude that $\Gamma(\R^2_{u,v})$ is either an open strip
between two parallel lines or an open half-plane. Here we use the
fact that invertible linear maps preserve half-planes and strips
bounded by parallel lines.

To express this classification intrinsically, suppose that a
Euclidean plane domain $D^2$ admits global Chebyshev coordinates
ranging over $\R^2$. The pullback of its Euclidean metric is an
entire flat Chebyshev metric. Its coordinate map and the map
$\Gamma$ are developing maps of the same flat metric, and hence
differ by a Euclidean rigid motion. The classification just proved
therefore applies to $D^2$.

\medskip
Next, we show the converse assertion. Let $h:\R\to(0,1)$ be smooth and set
\[
  \Omega_1(u):=\arcsin h(u),\qquad \Omega_2(v):=0.
\]
This is an $S$-pair, and its developing map is
\[
  \Gamma_h(u,v):=
  \left(v+\int_0^u\sqrt{1-h(s)^2}\,ds,
  -\int_0^u h(s)\,ds\right)^T.
\]
For each fixed $u$, its first component ranges over all of $\R$.
Its second component is strictly decreasing, so $\Gamma_h$ is a
diffeomorphism onto $\R\times I_h$, where
\[
  I_h:=\left\{-\int_0^u h(s)\,ds\,;\, u\in\R\right\}.
\]
Since $\Gamma_h^*(dx^2+dy^2)=d\sigma^2$, this diffeomorphism is
an isometry. The following choices give the three required types:
\begin{enumerate}
\item If $h(u)=1/2$, then $I_h=\R$.
\item If $h(u)=1/(2(1+e^u))$, then
\[
  -\int_0^u h(s)\,ds
  =\frac12\log\frac{1+e^{-u}}2,
  \qquad I_h=\left(-\frac{\log2}{2},\infty\right).
\]
\item Given $w>0$, if $h(u)=1/(2\cosh^2(u/w))$, then
\[
  -\int_0^u h(s)\,ds=-\frac w2\tanh(u/w),
  \qquad I_h=(-w/2,w/2).
\]
\end{enumerate}
Euclidean rigid motions give every half-plane and every strip of
width $w$, completing the converse.
\end{proof}

\begin{proof}[Proof of Theorem~B]
If $(\Omega_1,\Omega_2)$ is complete, Proposition~\ref{prop:821}
shows that $(\gamma_1,\gamma_2)$ is unbounded in all directions.
Conversely, suppose that $(\gamma_1,\gamma_2)$ is unbounded in all
directions.
By Theorem~\ref{thm:A}, $\Gamma(\R^2_{u,v})$ is either an open strip,
an open half-plane, or the entire plane.  If it were not the entire
plane, there would exist a unit vector $\mb e$ such that
$\Gamma\cdot\mb e$ is bounded above.  Since
\[
  \gamma_1(u)=\Gamma(u,0),\qquad
  \gamma_2(v)=\Gamma(0,v),
\]
both $\gamma_1\cdot\mb e$ and $\gamma_2\cdot\mb e$ would be bounded
above, contradicting Definition~\ref{def:all}.  Hence
$\Gamma(\R^2_{u,v})=\R^2_{x,y}$, and Proposition~\ref{prop:main}
shows that $(\Omega_1,\Omega_2)$ is complete.
\end{proof}

\begin{example}\label{ex:incomplete-over-another-plane}
A surface induced by an incomplete $S$-pair over one
space-like plane may be an entire graph over another space-like plane.
Set
\[
 \Omega_1(u):=0,
 \qquad
 \Omega_2(v):=\pi-\frac{\pi}{2(1+e^v)}.
\]
Then $\pi/2<\Omega_2(v)<\pi$, and hence
$(\Omega_1,\Omega_2)$ is an $S$-pair.  Put
\[
 X(v):=\int_0^v\cos\Omega_2(s)\,ds,
 \qquad
 Y(v):=\int_0^v\sin\Omega_2(s)\,ds.
\]
Its developing map is
\begin{equation}\label{eq:Gamma2}
 \Gamma(u,v)=\pmt{u+X(v)\\ Y(v)}.
\end{equation}
Since $\Omega_2(v)\to\pi/2$ as $v\to-\infty$, we have
$Y(v)\to-\infty$.  On the other hand,
\[
 0<\sin\Omega_2(v)
 =\sin\left(\frac{\pi}{2(1+e^v)}\right)
 \le \frac{\pi}{2}e^{-v}
 \qquad (v\ge0).
\]
Thus
\begin{equation}\label{eq:IA2}
 \Gamma(\R^2_{u,v})=\R\times(-\infty,B),
\qquad B:=\int_0^\infty\sin\Omega_2(v)\,dv<\infty.
\end{equation}
By Theorem~\ref{thm:A}, $(\Omega_1,\Omega_2)$ is not complete.

By \eqref{eq:564}
with \eqref{eq:Gamma2},
for $\sigma=1$, the corresponding surface is, up to a translation,
\[
 F(u,v)=\bigl(u+v,u+X(v),Y(v)\bigr).
\]
Fix $\epsilon\in(0,1)$ and consider the space-like plane
\[
 \Pi_\epsilon:=\bigl\{(z,x,y)\in\R^3_1\,;\, z+\epsilon y=0\bigr\}.
\]
Its Lorentz-normal direction is
$\mb n_\epsilon=(1,0,-\epsilon)$.  The functions
\[
 \xi:=x,
 \qquad
 \eta:=y+\epsilon(z-x)
\]
are constant along affine lines parallel to $\mb n_\epsilon$.  Moreover, on
$\Pi_\epsilon$ one has
\[
 \pmt{\xi\\ \eta}=\pmt{x\\ (1-\epsilon^2)y-\epsilon x},
\]
so $(\xi,\eta)^T$ is a global affine coordinate system on $\Pi_\epsilon$.
Hence the Lorentz-orthogonal projection of $F$ onto $\Pi_\epsilon$ is,
in these coordinates,
\[
 (u,v)\mapsto 
 \bigl(u+X(v),\Psi(v)\bigr),
 \qquad
 \Psi(v):=Y(v)+\epsilon\bigl(v-X(v)\bigr).
\]
Since $\cos\Omega_2(v)<0$, we have
\[
 \Psi'(v)
 =\sin\Omega_2(v)+\epsilon\bigl(1-\cos\Omega_2(v)\bigr)>\epsilon.
\]
Thus $\Psi$ is a diffeomorphism of $\R$ onto $\R$, and the projection
map above is a diffeomorphism, with inverse
\[
 v=\Psi^{-1}(\eta),
 \qquad
 u=\xi-X(v).
\]
Therefore $F$ is an entire graph over $\Pi_\epsilon$.  Since $\epsilon>0$ can be
chosen arbitrarily small, $\Pi_\epsilon$ can be chosen arbitrarily close to
the $xy$-plane.
\end{example}

\section{Applications}\label{sec:4}

\subsection*{Completeness of regular $S$-pairs}
We first give a proof of Fact~\ref{fact:227} using our
completeness criterion.

\begin{proposition}\label{prop:regular-complete}
If $(\Omega_1,\Omega_2)$ is a regular $S$-pair,
then it is complete.
\end{proposition}

\begin{proof}
Since $(\Omega_1,\Omega_2)$ is regular, there exists $\epsilon>0$ such that
\eqref{eq:Sin2} holds.  By \eqref{eq:733aaa}, there exist $\mu\in[0,\pi]$
and an $S$-pair $(\Omega_1^\star,\Omega_2^\star)$ which is $\R$-equivalent to
$(\Omega_1,\Omega_2)$ such that, for $\mu_0:=\mu+\epsilon/2$,
\begin{equation}\label{eq:733ax2}
\epsilon \le \hat \Omega_1^\star \le \mu_0-\frac{\epsilon}{2}
<\mu_0<\mu_0+\frac{\epsilon}{2}\le \Omega_2^\star \le \pi .
\end{equation}
In particular, $0<\mu_0<\pi$.
We denote by $\gamma_1^\star$ and $\gamma_2^\star$ the developing curves associated with
the $S$-pair $(\Omega_1^\star,\Omega_2^\star)$.
Then
\[
\bar\gamma_i=\mc L_{\mu_0}\circ\gamma_i^\star
\qquad (i=1,2)
\]
defines a monotone pair in the coordinates associated with
$\mu_0$. Since
$\epsilon \le \hat \Omega_1^\star \le \mu_0-\epsilon/2$,
the curve $\bar\gamma_1$ is the graph of an entire function
over the $\xi$-axis and is unbounded both above and below.
Corollary~\ref{cor:1436} and Proposition~\ref{prop:714}
therefore show that $(\Omega_1,\Omega_2)$ is complete.
\end{proof}

\subsection*{Associated families of time-like entire ZMC graphs}

\begin{proof}[Proof of Corollary~C]
Since $u\mapsto e^tu$ and $v\mapsto e^{-t}v$ are bijections of $\R$,
the pair
\[
 (\Omega_1(e^tu),\Omega_2(e^{-t}v))
\]
is an $S$-pair whenever $(\Omega_1,\Omega_2)$ is.  Its developing
curves are
\[
 \tilde\gamma_1(u)=e^{-t}\gamma_1(e^tu),
 \qquad
 \tilde\gamma_2(v)=e^t\gamma_2(e^{-t}v).
\]
Each curve is obtained by a reparametrization followed by multiplication
by a positive constant.  These operations preserve properness and the
property that a given scalar projection is unbounded above.  Hence
$(\tilde\gamma_1,\tilde\gamma_2)$ is unbounded in all directions
if and only if $(\gamma_1,\gamma_2)$ is.  The completeness assertion now
follows from Theorem~B.
\end{proof}

\subsection*{Conjugate surfaces need not be entire}

Up to a translation, write $F(u,v)=\alpha(u)+\beta(v)$ as in
\eqref{eq:564}.  By the definition of the conjugate surface, it is
represented by $\alpha(u)-\beta(v)$.  After the reparametrization
$v\mapsto-v$, the same surface is represented by
\[
 F^*(u,v)=\alpha(u)-\beta(-v),
\]
and hence
\[
 (F^*)_u
 =\bigl(\sigma,\cos\Omega_1(u),-\sin\Omega_1(u)\bigr),
 \qquad
 (F^*)_v
 =\bigl(\sigma,\cos\Omega_2(-v),\sin\Omega_2(-v)\bigr).
\]
Thus $F^*$ is induced by the $S$-pair
$(\Omega_1(u),\Omega_2(-v))$.

\begin{proof}[Proof of Theorem~D]
Since $v\mapsto-v$ is a bijection of $\R$, the conjugate pair is
regular if and only if $(\Omega_1,\Omega_2)$ is.  Its developing map
is the spatial projection of $F^*$.  Theorem~\ref{thm:A} therefore
shows that $F^*$ is a graph over an open strip, an open half-plane, or
the entire $xy$-plane.  If the original pair is regular, the conjugate
pair is complete by Proposition~\ref{prop:regular-complete}, and
$F^*$ is entire by Proposition~\ref{cor:374}.  The final assertion
follows from Proposition~\ref{prop:conjugate-no-spacelike-graph}
below.
\end{proof}

Set
\begin{equation}\label{eq:1876}
 L(t):=\log(1+e^{-t}),\qquad
 \Omega(t):=\frac{L(t)}{1+L(t)^2}.
\end{equation}

\begin{proposition}\label{prop:conjugate-no-spacelike-graph}
The pair $(\Omega,\Omega)$ is a complete but non-regular $S$-pair,
and the corresponding time-like ZMC surface is an entire graph.
However, its conjugate surface is not an entire graph over any
space-like affine plane in $\R^3_1$.
\end{proposition}

\begin{figure}[h!]
  \begin{center}
    \includegraphics[width=\linewidth,height=56mm,keepaspectratio]{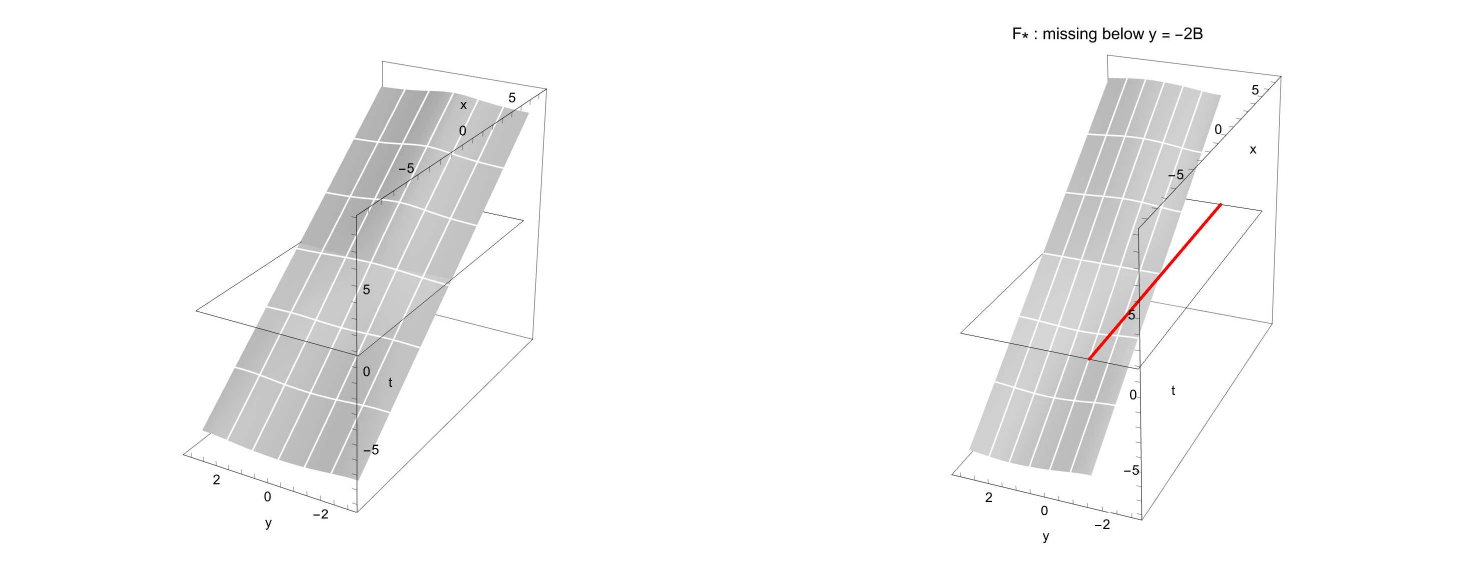}
    \caption{The time-like entire ZMC graph $F$ and its conjugate surface $F^*$. 
The surface $F^*$ terminates along a light-like line at $y=-2B\approx -1.39984$ 
and is not an entire graph.}
    \label{fig2}
  \end{center}
\end{figure}

\begin{proof}
We have $0<\Omega\leq1/2$ and
\[
 \Omega(t)\sim e^{-t}\quad(t\to+\infty),
 \qquad
 \Omega(t)\sim |t|^{-1}\quad(t\to-\infty),
\]
where $\Omega(t)\sim e^{-t}$ 
(resp. $\Omega(t)\sim |t|^{-1}$) 
as $t\to+\infty$ (resp. $t\to-\infty$)
means that
\[
\lim_{t\to+\infty}\frac{\Omega(t)}{e^{-t}}=1
\qquad \left(
\text{resp. }\, \lim_{t\to-\infty}\frac{\Omega(t)}{|t|^{-1}}=1\right).
\]
Thus $0<\Omega(u)+\Omega(v)\leq1<\pi$, so $(\Omega,\Omega)$
is an $S$-pair.  It is not regular since
$\sin(2\Omega(t))\to0$ as $t\to+\infty$.

Set
\begin{equation}\label{eq:125}
 X(t):=\int_0^t\cos\Omega(s)\,ds,
 \qquad
 Y(t):=\int_0^t\sin\Omega(s)\,ds.
\end{equation}
Since $\sin\Omega(s)\sim e^{-s}$ as $s\to+\infty$, we may define
\[
 B:=\int_0^\infty\sin\Omega(s)\,ds\in(0,\infty).
\]
Then $Y(t)\to B$ as $t\to+\infty$, and
\[
 B-Y(t)=\int_t^\infty\sin\Omega(s)\,ds>0
 \qquad(t\in\R).
\]

The developing curves are
\[
 \gamma_1(t)=\pmt{X(t)\\ -Y(t)},
 \qquad
 \gamma_2(t)=\pmt{X(t)\\ Y(t)}.
\]
Since $X'\geq\cos(1/2)>0$, both curves are simple and proper,
and $X(t)\to\pm\infty$ as $t\to\pm\infty$.
Moreover, $\sin\Omega(t)\sim |t|^{-1}$ as $t\to-\infty$,
so $Y(t)\to-\infty$ as $t\to-\infty$.
Consequently, for every unit vector $\mb e=(a,b)^T$,
\[
 \max\{\gamma_1(t)\cdot\mb e,\gamma_2(t)\cdot\mb e\}
 =aX(t)+|bY(t)|
\]
is unbounded above.  Hence $(\gamma_1,\gamma_2)$ is unbounded
in all directions.  By Theorem~B and
Proposition~\ref{cor:374}, the pair $(\Omega,\Omega)$ is complete
and the corresponding surface is an entire graph over the $xy$-plane.

Up to a Lorentzian rigid motion, the conjugate surface is
\[
 F^*(u,v)
 =\bigl(u+v,\,X(u)-X(-v),\,-Y(u)-Y(-v)\bigr).
\]
Since $Y(t)<B$ for all $t\in\R$, its image lies in
$\{y>-2B\}$.  On the other hand,
\[
 F^*(s,-s)=(0,0,-2Y(s))
 \rightarrow (0,0,-2B)
 \qquad(s\to+\infty).
\]
Thus $F^*(\R^2)$ is not closed in $\R^3_1$.
Since an entire graph of a continuous function over any space-like 
affine plane is closed, the conjugate surface cannot be such a graph 
(see Figure \ref{fig2}).
\end{proof}

\subsection*{Entire ZMC graphs with unbounded Gaussian curvature}

We first construct an entire ZMC graph with unbounded Gaussian
curvature arising from a regular $S$-pair.

\begin{example}\label{ex:regular-unbounded-K}
The Gaussian curvature of a time-like entire ZMC graph of type $S$
can be unbounded even when its associated $S$-pair is regular.
Indeed, set
$$
 \Omega(t):=\frac{\pi}{4}+\frac12\sin(t^2),
 \qquad
 \Omega_1(u):=\Omega(u),\qquad
 \Omega_2(v):=\Omega(v).
$$
Then
$$
 \frac{\pi}{2}-1
 \leq \Omega_1(u)+\Omega_2(v)
 \leq \frac{\pi}{2}+1
 \qquad (u,v\in\R).
$$
Since $1<\pi/2$, it follows that
$$
 \sin\bigl(\Omega_1(u)+\Omega_2(v)\bigr)
 \geq\cos1>0.
$$
Thus $(\Omega_1,\Omega_2)$ is a regular $S$-pair.  By
Proposition~\ref{prop:regular-complete} and
Proposition~\ref{cor:374}, the corresponding ZMC immersion is a
time-like entire ZMC graph of type $S$.

On the other hand,
$
 \Omega'(t)=t\cos(t^2).
$
Put $t_n:=\sqrt{n\pi}$ for $n\geq1$.  Then
$$
 \sin(t_n^2)=0,
 \qquad
 \cos(t_n^2)=(-1)^n,
$$
and hence
$
 \Omega_1(t_n)+\Omega_2(\pm t_n)=\pi/2.
$
Using \eqref{eq:K-Magid}, we obtain
$$
 K(t_n,t_n)=-t_n^2=-n\pi,
 \qquad
 K(t_n,-t_n)=t_n^2=n\pi.
$$
Consequently, the Gaussian curvature is unbounded both above and
below.
\end{example}

In this example, both $\Omega_1'$ and $\Omega_2'$ are unbounded.
We now prove that the same curvature behavior can occur even when
both derivatives are bounded.

\begin{proof}[Proof of Theorem~F]
Let $F$ be the time-like entire ZMC graph associated with the
complete but non-regular $S$-pair $(\Omega,\Omega)$ in
Proposition~\ref{prop:conjugate-no-spacelike-graph}, where $\Omega$
is defined in \eqref{eq:1876}.  We show that its Gaussian curvature
is unbounded both above and below.

By \eqref{eq:1876},
$$
 L'(t)=-\frac{1}{1+e^t},\qquad
 \Omega'(t)
 =
 \frac{L(t)^2-1}
 {(1+e^t)(1+L(t)^2)^2}.
$$
In particular,
$$
 \Omega(t)\sim e^{-t},
 \qquad
 \Omega'(t)\sim-e^{-t}
 \qquad (t\to+\infty),
$$
whereas
$$
 \Omega(t)\sim\frac1{|t|},
 \qquad
 \Omega'(t)\sim\frac1{t^2}
 \qquad (t\to-\infty).
$$
Since $\Omega'$ is smooth and tends to zero at both ends of $\R$, it
is bounded.

By \eqref{eq:K-Magid},
$$
 K_F(u,v)=
 -\frac{\Omega'(u)\Omega'(v)}
 {\bigl(1-\cos(\Omega(u)+\Omega(v))\bigr)^2}.
$$
Using
the fact
$
 (1-\cos x)\sim x^2/2
$
as $x\to0$, we obtain, along $u=v=t$,
$$
 K_F(t,t)
 \sim
-\frac14e^{2t}
 \rightarrow-\infty
 \qquad (t\to+\infty).
$$
On the other hand,
$$
 \Omega(-e^t)\sim e^{-t},
 \qquad
 \Omega'(-e^t)\sim e^{-2t}
 \qquad (t\to+\infty).
$$
Therefore,
$$
 K_F(t,-e^t)
 \sim
\frac14e^t
 \rightarrow+\infty.
$$
Thus the Gaussian curvature of $F$ is unbounded both above and below.
Since $\Omega_1'=\Omega_2'=\Omega'$ is bounded, the theorem follows.
\end{proof}

\begin{acknowledgements}
The authors thank Professor Yoshihisa Kitagawa for valuable comments.
\end{acknowledgements}

\section*{Use of generative AI}
GPT-6-Astra (OpenAI) was used to improve the language and clarity of 
the manuscript. The authors take full responsibility for the manuscript.


\begin{thebibliography}{99}

\bibitem{Ac}
S.~Akamine,
{\it Isometric and anti-isometric classes of timelike minimal surfaces
in Lorentz--Minkowski space}, Differential Geom. Appl. {\bf 97}
(2024), Paper No. 102210.

\bibitem{AF}
S.~Akamine and H.~Fujino,
\newblock \emph{Extension of Krust theorem and deformations of minimal
surfaces},
\newblock Ann. Mat. Pura Appl. (4)
\textbf{201} (2022), 2583--2601.

\bibitem{AGM}
J.~A.~Aledo, J.~A.~G\'alvez and P.~Mira,
{\it A D'Alembert formula for flat surfaces in the $3$-sphere},
J. Geom. Anal. {\bf 19} (2009), 211--232.

\bibitem{Karcher}
H.~Karcher,
\newblock \emph{Construction of minimal surfaces},
\newblock in \emph{Surveys in Geometry}, University of Tokyo,
pp.~1--96, 1989.
\newblock Also published as Lecture Notes No.~12, SFB~256,
University of Bonn, 1989.

\bibitem{Ki}
Y.~Kitagawa, {\it
Periodicity of the asymptotic curves on
flat tori in $S^3$}, J. Math. Soc.\ Japan {\bf 40} (1988), 457--476.

\bibitem{Kon}
J.~J.~Konderak,
{\it A Weierstrass representation theorem for Lorentz surfaces},
Complex Var.\ Theory Appl. {\bf 50} (2005),  319--332.

\bibitem{LMP}
M.~A.~Le\'on-Guzm\'an, P.~Mira
and J.~A.~Pastor,
{\it The space of Lorentzian flat tori
in anti-de Sitter 3-space},
Trans.\ Amer.\ Math.\ Soc.\
{\bf 363} (2011), 6549--6573.

\bibitem{Lopez}
R.~L\'opez,
\newblock \emph{The Lorentzian version of a theorem of Krust},
\newblock Bull. Malays. Math. Sci. Soc.
\textbf{44} (2021), 683--692.

\bibitem{M}
M.~A.~Magid,
{\it The Bernstein problem for timelike surfaces},
Yokohama Math.\ J.\ {\bf 37}
(1989), 125--137.

\bibitem{W0}
T.~K.~Milnor,
{\it Entire timelike minimal surfaces in $E^{3,1}$},
Michigan Math.\ J.\ {\bf 37}
(1990), 163--177.

\bibitem{W}
T.~K.~Milnor,
{\it Associate harmonic immersions in $3$-space},
Michigan Math.\ J.\ {\bf 37}
(1990), 179--190.

\end{thebibliography}
\end{document}